\documentclass{elsart}
\usepackage{amssymb}
\usepackage{stackrel}
\usepackage{amsmath, mathtools}
\usepackage{amsfonts}
\usepackage{graphicx}

\begin{document}

\begin{frontmatter}



{\title{Size of a $3$-uniform linear hypergraph}}
\begin{abstract}
This article provides bounds on the size of a $3$-uniform linear hypergraph with restricted matching number and maximum degree. In particular, we show that if a $3$-uniform, linear family $\mathcal{F}$ has maximum matching size $\nu$ and maximum degree $\Delta$ such that $\Delta\geq \frac{23}{6}\nu(1+\frac{1}{\nu-1})$, then $|\mathcal{F}|\leq \Delta \nu$. 
\end{abstract}
\author{Niraj Khare}\footnote{This article is dedicate to Prof. Akos Seress for his contribution to Combinatorics and inspiring love for the field in his students.}
\ead{nirajkhare@math.ohio-state.edu }
\address{The Ohio State University, Columbus, Ohio, USA }
\end{frontmatter}
\section{Introduction}
Let $V$ be a set of vertices and let $\mathcal{F}\subseteq 2^{V}$ be a set of distinct subsets of $V$. A set system $\mathcal{F}$ is $k$-uniform for a positive integer $k$ if $|A|=k$ for all $A \in \mathcal{F}$. A set system $\mathcal{F}$ is linear if $|A \cap B|\leq 1$ for all distinct $A$, $B$ in $\mathcal{F}$. For a hypergraph $\mathcal{G}=(V, \mathcal{F})$, the set $V$ is called the set of vertices of $\mathcal{G}$ and the set $\mathcal{F}\subseteq 2^{V}$ is called the set of hyper-edges of $\mathcal{G}$. The size of a $k$-uniform linear hypergraph $\mathcal{G}=(V, \mathcal{F})$ is $|\mathcal{F}|$-the number of its hyper edges. A matching in $\mathcal{G}$ (or $\mathcal{F}$) is a collection of pairwise disjoint hyper-edges of $\mathcal{G}$. The size of a maximum matching in $\mathcal{F}$ shall be denoted by $\nu(\mathcal{F})$. Also degree of a vertex and maximum degree of $\mathcal{G}$ is defined in a usual familiar way. For any $x \in V$, define $\mathcal{F}_x = \{ A \in \mathcal{F}\quad|\quad x \in A\}$ and $\Delta(\mathcal{F})=\max\{|\mathcal{F}_x|\quad|\quad x \in V\}$. The objective of this article is to find a bound on the size of $\mathcal{F}$ for given values of $\Delta(\mathcal{F})$ and $\nu(\mathcal{F})$. Throughout the remainder of this article unless otherwise stated, $\mathcal{F}$ shall be a $3$-uniform linear set system with maximum matching size $\nu(\mathcal{F})=\nu$ and maximum degree $\Delta(\mathcal{F})=\Delta$. Also, for any set system $\mathcal{H}$ and $\mathcal{B}\subseteq \mathcal{H}$, we shall use the following notation: $X_{\mathcal{B}}:= \stackrel [A \in \mathcal{B}]{}{\cup}A $.

The problem of bounding the size of a uniform family by restricting matching size and maximum degree has been studied for simple graphs in \cite{CH} and \cite{bk}. These articles were in turn inspired by the \textit{Sunflower} lemma due to Erd\H{o}s and Rado (see \cite{ER}). A \textit{sunflower} with $s$ petals is a collection of sets $A_1,A_2\ldots, A_s$ and a set $X\textrm{(possibly empty)}$ such that $A_i\cap A_j=X$ whenever $i\neq j$. The set $X$ is called the core of the sunflower. A linear family admits two kinds of \textit{Sunflower}: (i) a matching is a \textit{Sunflower} with an empty core, (ii) a collection of hyper-edges incident at a vertex. It is a well known result (due to Erd\H{o}s-Rado\cite{ER}) that a $k$-uniform set system, with more members than $k!(s-1)^k$ admits a sunflower with $s$ petals (for a proof see \cite{BF}). Other bounds that ensure the existence of a sunflower with $s$ petals are known in the case of $s=3$ with block size $k$ (see \cite{Ka}). However, not much progress has been made towards the general case. This article considers the dual problem of finding the maximum size of a $3$-uniform, linear family $\mathcal{F}$ that admits no \textit{Sunflower} with $s$ petals, i.e., $s> \nu(\mathcal{F})$ and $s> \Delta(\mathcal{F})$. The following remark on the size of a family shall be useful later.  
\begin{rem}\label{delta-nu}
For a positive integer $\Delta$, let a $3$-uniform family $\mathcal{G}$ be a \textit{Sunflower} with $\Delta$ petals and core of size one. For any positive integer $\nu$, let $\mathcal{F}$ consists of $\nu$ components where each component is isomorphic to $\mathcal{G}$. It is obvious that $\nu(\mathcal{F})=\nu$ and $\Delta(\mathcal{F})=\Delta$. Also, $|\mathcal{F}|=\Delta \nu$.
\end{rem}
\section{Results}
Our aim in this article is to prove the following two results.
\begin{thm}\label{thm0} 
Let $\mathcal{F}$ be a $3$-uniform linear set system with maximum matching size $\nu(\mathcal{F})=\nu$ and maximum degree $\Delta(\mathcal{F})=\Delta$. If $\Delta\geq 5$, then $|\mathcal{F}|\leq 2\Delta \nu$.
\end{thm}
The main result of this article is a tighter bound in case $\Delta$ is approximately greater than $4\nu$. The precise statement follows. 
\begin{thm} {(\bf{The main result})}\label{thm1}
Let $\mathcal{F}$ be a $3$-uniform linear set system with maximum matching size $\nu(\mathcal{F})=\nu$ and maximum degree $\Delta(\mathcal{F})=\Delta$. If $\Delta\geq \frac{23}{6}\nu(1+\frac{1}{\nu-1})$, then $|\mathcal{F}|\leq \Delta \nu$.\end{thm}
Let $\nu$ be any positive integer. It is worthwhile to note that there are $3$-uniform liner families $\mathcal{F}$ with $\nu=\nu(\mathcal{F})$ such that $|\mathcal{F}|>\Delta(\mathcal{F})\nu(\mathcal{F})$. In the next section we construct such families and thus establish importance of the main result-Theorem \ref{thm1}.
\section{Families with large size}  
Let $\mathcal{F}$ be a $3$-uniform linear family with $\Delta:=\Delta(\mathcal{F})$ and $\nu:=\nu(\mathcal{F})$. We present some examples such that $|\mathcal{F}|> \Delta \nu$.\newline
(i) There are block designs $\mathcal{F}$ with block size three such that $|\mathcal{F}|\geq\nu(\mathcal{F})\Delta(\mathcal{F})$. For example, consider Steiner triples $S(n,3,2)$. A Steiner system $S(n,k,r)$ is a set system on $n$ vertices such that each member has cardinality $k$ and every $r$-subset of vertices is contained in a unique member (also called block) of the family $S(n,k,r)$. It is well known that $S(n,3,2)$ exists if and only if $n\geq 3$, and $n\equiv 1\mbox{ mod($6$)}$ or $n\equiv 3\mbox{ mod($6$)}$ (see \cite{CR}, for instance). 
\begin{itemize}
\item If $n=6m+1$ and $\mathcal{F}$ is an $S(n,3,2)$ then $|\mathcal{F}|=\frac{1}{3}{6m+1 \choose 2}=m(6m+1)$, $\Delta(\mathcal{F})=3m$, and $\nu(\mathcal{F})\leq 2m$, so $|\mathcal{F}|>\Delta(\mathcal{F})\nu(\mathcal{F})$.
\end{itemize}
(ii) By the method given in \cite{bk}, we can construct a simple graph $G$ for any $\Delta:=\Delta(G)$ and $\nu:=\nu(G)$ such that $|E(G)|=\nu\Delta+\lfloor{\frac{\nu}{\lceil{\frac{\Delta}{2}}\rceil}}\rfloor \lfloor{\frac{\Delta}{2}}\rfloor$. Note that if $2\leq \Delta \leq 2\nu$ then $|E(G)|>\Delta \nu$. Let $Y$ be a set such that $Y\cap V(G)=\emptyset$ and $|Y|=|E(G)|$. We order the edges $\{e_1,e_2,\ldots,e_{|E(G)|}\}$ in $E(G)$ randomly and let $Y=\{y_1,y_2,\cdots,y_{|E(G)|}\}$. We define a linear, $3$-uniform family $\mathcal{F}$ such that $\nu(\mathcal{F})=\nu(G)$ and $\Delta(\mathcal{F})=\Delta(G)$. For $i \in \{1,2,\ldots,|E(G)|\}$, let $A_i:=e_i\cup \{y_i\}$. Now let $\mathcal{F}:=\{A_i|\ i\in \{1,2,\ldots,|E(G)|\}\}$. It is obvious that $\mathcal{F}$ is a $3$-uniform, linear family. Also note that $\nu(\mathcal{F})=\nu$, $\Delta(\mathcal{F})=\Delta$ and $|\mathcal{F}|=|E(G)|$. Thus, $|\mathcal{F}|=|E(G)|=\nu\Delta+\lfloor{\frac{\nu}{\lceil{\frac{\Delta}{2}}\rceil}}\rfloor \lfloor{\frac{\Delta}{2}}\rfloor > \Delta \nu$.

Theorem \ref{thm1} states that if $\Delta$ is large enough compared to $\nu$ then $|\mathcal{F}|\leq \nu\Delta$. On the other hand the example in part (ii) above shows that for any positive integer $\nu$, there are families $\mathcal{F}$ such that $|\mathcal{F}|>\Delta \nu$ with $2\leq\Delta \leq 2\nu$. It would be interesting to determine the exact value $f(\nu)$ so that for any $3$-uniform , linear family $\mathcal{F}$ with $\Delta(\mathcal{F})=\Delta \geq f(\nu)$ and $\nu(\mathcal{F})=\nu$, we have $|\mathcal{F}|\leq \nu \Delta$. 
\section{Preliminaries}
We first find a trivial bound to establish that the problem is well founded. Let $\mathcal{H}$ be a $k$-uniform set system with maximum matching $\nu$ and maximum degree $\Delta$. Since the set of vertices that are covered by a maximum matching form a vertex cover (also known as transversal), each hyperedge is covered by $k \nu$ vertices. As the maximum degree is $\Delta$, we get \begin{equation} |\mathcal{H}|\leq (\Delta-1)(k \nu)+\nu.  \label{trivial-bound}\end{equation} In general this bound is too large and can be improved. Surprisingly for $k=3$, there are values of $\nu$ and $\Delta$ for which the previous crude bound is tight. For example Fano plane of order two achieves the bound for $k=3$, $\Delta=3$ and $\nu=1$. Note that for $\Delta=2$ and $k=3$, the set system $\{ \{x,y,z\}, \{a,c,z\}, \{a, b,x\}, \{b,c,y\} \}$ on vertices $\{x,y,z,a,b,c\}$ satisfies eq (\ref{trivial-bound}).  Our aim is to improve the bound in eq (\ref{trivial-bound}) to obtain results of Theorem \ref{thm0} and Theorem \ref{thm1}. One of the critical lemmas needed is Lemma \ref{augmenting set lem}. This lemma is a generalized version of augmenting path maximum matching lemma for graphs. The statement of the augmenting path maximum matching lemma for graphs is that a matching is maximum if and only if there is no augmenting path relative to it. Readers can find graph theoretic version in any standard text book such as \cite{DW} or \cite{LP}. There are numerous versions available that extend augmenting path maximum matching lemma to hypergraphs (see \cite{F},for instance). However, the version presented here (i.e., Lemma \ref{augmenting set lem}) suits to our requirements better. Note that Lemma \ref{augmenting set lem} holds for any hypergraph and we don't require uniformity of cardinality of hyperedges. 
\begin{defn}\label{augmenting set}
\textit {Augmenting set} : Let $\mathcal {F}$ be a set system with a matching $\mathcal {M}$. We say $\mathcal {C} \subseteq \mathcal {F}$ is an $\mathcal{M}$-\textit {augmenting set} if and only if  $\mathcal {C}$ satisfies:\\
$[(1)]$ $|\mathcal {M} \cap \mathcal {C}| < |\mathcal {C} \setminus \mathcal {M}|$,\newline
       (i.e., there are more non-matching edges than matching edges in $\mathcal {C}$)\\ 
$[(2)]$ If $B \in \mathcal {M} $, $B \cap A \neq \emptyset $ for some $A \in \mathcal{C}$ then $B \in \mathcal {C}$,\newline
(i.e., if any matching edge has a non-empty intersection with any of the non-matching edges of $\mathcal {C}$ than that matching edge is also in $\mathcal {C}$)\\
$[(3)]$ $| \mathcal {C}_x \setminus \mathcal{M}| \leq 1$   $\forall x \in X_{\mathcal {C}}=\stackbin[A \in \mathcal{C}]{}{\cup}A $.\newline
    (i.e., any vertex of $\mathcal {C}$ is covered by at most one  \underline {non-matching edge of $\mathcal {C}$} or in other words, non-matching edges in $\mathcal{C}$ are pairwise disjoint.)
\end{defn}    
\begin{lem}\label{augmenting set lem}    
 Let $\mathcal {F}$ be a hypergraph and $\mathcal {M}$ be a matching. $\mathcal {M}$ is maximum if and only if there is no $\mathcal {M}$-\textit {augmenting set} in $\mathcal {F}$.
\end{lem} 
\noindent{\sl Proof.} We first show the only if part by proving the contrapositive. Suppose there is an $\mathcal {M}$ \textit {augmenting set} $\mathcal {C}$ in $\mathcal {F}$. Then we define a new subfamily, $\mathcal {M}_1 := \{\mathcal {M} \setminus \mathcal {C}\}  \cup \{\mathcal {C} \setminus \mathcal {M}\}$. Note that $| \mathcal {M}_1|> |\mathcal{M}|$ as $|\mathcal {C} \setminus \mathcal {M}| > |\mathcal{C} \cap \mathcal {M}|$ by property (1) of \textit{augmenting set}, Definition \ref{augmenting set}. We claim $\mathcal {M}_1$ is a matching of $\mathcal{F}$. Note that two non-matching edges of $\mathcal {C}$ do not intersect by the property (3) of \textit {augmenting set} (Definition \ref{augmenting set}), and no edge of $\mathcal {M} \setminus \mathcal {C}$ can have non-empty intersection with an edge of $\mathcal {C}$ by the property (2) of  \textit {augmenting set} (Definition \ref{augmenting set}). Also edges in $\mathcal {M} \setminus \mathcal {C}$ are pairwise disjoint as $\mathcal{M}$ is a matching. Therefore, members of $\mathcal {M}_1$ are pairwise disjoint. Thus, $\mathcal {M}_1$ is a matching of $\mathcal{F}$.\newline
Next we prove the if part. Let $\mathcal {M}$ be a matching of $\mathcal {F}$ which is not maximum and $\mathcal {M}_1$ be a maximum matching. Hence $|\mathcal {M}_1|> |\mathcal {M}|$. Let $\mathcal {S} := \{\mathcal {M}_1 \setminus \mathcal{M} \} \cup \{\mathcal {M} \setminus \mathcal{M} _1\}$. In $\mathcal{S}$ there are more $\mathcal{M}_1$ edges than $\mathcal{M}$ edges. So there exists a component $\mathcal {C}$ of $\mathcal {S}$ such that $\mathcal{C}$ contains more $\mathcal{M}_1$ edges than $\mathcal{M}$ edges. We claim that $\mathcal {C}$ is an $\mathcal{M}$-\textit {augmenting set} by Definition \ref{augmenting set} as,\newline
[(1)] $\mathcal{C}$ has more non-matching (relative to $\mathcal {M}$) edges than matching $\mathcal {M} $ edges,\newline
[(2)] $\mathcal {C}$ is a component hence any $\mathcal{M}$ edge which has a non-empty intersection with any of the $\mathcal {C}$ edges is in $\mathcal {C}$. Note that any edge in $\mathcal{M}_1 \cap \mathcal {M}$ can not have non-empty intersection with any of the $\mathcal {C}$ edges,\newline
[(3)] $ | \mathcal {C}_x \setminus \mathcal{M}| \leq 1$   $\forall x \in X_{\mathcal {C}}$ holds trivially as $\mathcal {M}_1$ is a matching of $\mathcal{F}$. \hspace{\stretch{1}} $\square$\newline

It is easy to prove the first result, i.e., Theorem \ref{thm0}. However, some more definitions are needed to this end. Let $\mathcal{M}$ be a maximum matching of a $k$-uniform family $\mathcal{F}$. For $i\in\{1,2,\cdots,k\}$, define
$D_i(\mathcal{F, \mathcal{M}}):= \{ A\in \mathcal{F}\quad|\quad |A\cap X_{\mathcal{M}}|=i\}$.
Also we define for $x\in X_{\mathcal{F}}$,
$d_i(x, \mathcal{M}):= | \{A\in D_i(\mathcal{F}, \mathcal{M})\quad|\quad x\in A\}|$ for $i\in \{1,2,\cdots,k\}$.\newline
\begin{lem}\label{degree 0-1} Let $\mathcal{F}$ be a linear $k$-uniform family with $k\geq 2$ and $\mathcal{M}$ be a maximum matching of $\mathcal{F}$. If $B=\{x_1,x_2,\ldots,x_k\}$ is an $\mathcal{M}$ edge such that for some $1\leq i\leq k$, $d_1(x_i, \mathcal{M})\geq k$ then $d_1(x_j, \mathcal{M})=0$ for all $j\neq i$ and $1\leq j\leq k$.
\end{lem}
\noindent{\sl Proof.} Without loss of generality, let $i=1$ and let $\mathcal{F}_{x_1}\cap D_1(\mathcal{F}, \mathcal{M})=\{A_i\quad|\quad i\in I\}$ where $|I|=d_1(x_1,\mathcal{M})\geq k$. As $\mathcal{F}$ is a linear family, we have $\cap_{i\in I}A_i=\{x_1\}$. Suppose on the contrary $d_1(x_j, \mathcal{M})\geq 1$ for some $j\neq 1$. Let $C\in D_1(\mathcal{F})\cap \mathcal{F}_{x_j}$. As $|I|\geq k$, the sets $A_i\setminus \{x_1\}$ are pairwise disjoint for $i\in I$ and $|C\setminus \{x_j\}|=k-1$, linearity of $\mathcal{F}$ demands that $C \cap A_i=\emptyset$ for some $i\in I$.  By Definition \ref{augmenting set}, $\{C,A_i,B\}$ is an $\mathcal{M}$-augmenting set since the only matching edge covered by $A_i$ and $C$ is $B$ and $C\cap A_i=\emptyset$. It is a contradiction to Lemma \ref{augmenting set lem} as $\mathcal{M}$ is a maximum matching. \hspace{\stretch{1}} $\square$\newline
\begin{lem}\label{D1 edges} Let $k\geq 2$ be a positive integer. If $\mathcal{F}$ be a linear $k$-uniform family with a maximum matching $\mathcal{M}$ then\newline 
 $|D_1(\mathcal{F},\mathcal{M})|\leq \max\{(\Delta-1)\nu, k(k-1)\nu\}$.
\end{lem}
\noindent{\sl Proof.} For $B\in \mathcal{M}$, let $\mathcal{D}_1(B):=\{A\in D_1(\mathcal{F},\mathcal{M})\quad|\quad A\cap B \neq \emptyset\}$. It is enough to show that for each $B\in \mathcal{M}$, $|\mathcal{D}_1(B)|\leq \max\{\Delta-1,k(k-1)\}$.\newline
Suppose that for $B=\{x_1,\ldots,x_k\}\in \mathcal{M}$, $|\mathcal{D}_1(B)|\geq k(k-1)+1$. Then there exists, by pigeon hole principal, a $x_i \in B$ contained in at least $k$ members of $\mathcal{D}_1(B)$. Thus, by Lemma \ref{degree 0-1} all $D_1(B)$ edges are incident at $x_i$ (i.e.,  $X_{\mathcal{D}_1(B)}\cap B = \{x_i\}$). Since $x_i$ is contained in at most $\Delta-1$ elements of $\mathcal{F}$ different from $B$, we obtain $|\mathcal{D}_1(B)|\leq \Delta-1$. \hspace{\stretch{1}} $\square$\newline 
Next we rewrite and prove Theorem \ref{thm0} using the last lemma. 
\begin{thm}\label{2dnu bound}
Let $\mathcal{F}$ be a linear $3$-uniform family. If $\Delta(\mathcal{F})=d$ and $\nu(\mathcal{F})=\nu$ then
\begin{equation}
|\mathcal{F}|\leq \max\{2d\nu, 10\nu\}.
\end{equation}
\end{thm}
\noindent{\sl Proof.} Let $\mathcal{M}$ be a maximum matching of $\mathcal{F}$. For any $k$-uniform family, the summation of degrees of vertices is equal to $k$ times the number of edges. Hence for $k=3$,
\begin{equation}\label{degree-edge-count}
\stackrel[{x \in (X_{\mathcal{F}}\setminus X_{\mathcal{M}})}]{}{\sum}|\mathcal{F}_x|+ \stackrel[{x \in X_{\mathcal{M}}}]{}{\sum}|\mathcal{F}_x|=3|\mathcal{F}|.
\end{equation}
Now we consider the following two cases.\newline
{\bf{Case I}}: $\stackrel[{x \in (X_{\mathcal{F}}\setminus X_{\mathcal{M}})}]{}{\sum}|\mathcal{F}_x|\leq \stackrel[{x \in X_{\mathcal{M}}}]{}{\sum}|\mathcal{F}_x|$.\newline
By equation (\ref{degree-edge-count}) and the case assumption, 
\begin{displaymath}2\stackrel[{x \in X_{\mathcal{M}}}]{}{\sum}|\mathcal{F}_x|\geq 3|\mathcal{F}|. \end{displaymath}
As $|X_{\mathcal{M}}|=3\nu$, we have $\stackrel[{x \in X_{\mathcal{M}}}]{}{\sum}|\mathcal{F}_x|\leq d|X_{\mathcal{M}}|=3d\nu$. Therefore,\newline
\begin{displaymath}
2 (3 d \nu)\geq 2{\sum}_{x \in X_{\mathcal{M}}}|\mathcal{F}_x|\geq 3|\mathcal{F}|.
\end{displaymath} Thus,
\begin{equation}\label{case_1_eq}
 2d\nu\geq |\mathcal{F}|.
\end{equation}
{\bf{Case II}}: $\stackrel [{x \in (X_{\mathcal{F}}\setminus X_{\mathcal{M}})}]{}{\sum}|\mathcal{F}_x|>\stackrel[{x \in X_{\mathcal{M}}}]{}{\sum}|\mathcal{F}_x|$.\newline
As before, for $i \in \{1, 2,3\}$ define $D_i(\mathcal{F},\mathcal{M}):=\{A\in \mathcal{F}\quad|\quad |A \cap X_{\mathcal{M}}|=i\}$  and  $d_i:=|D_i(\mathcal{F},\mathcal{M})|$. Note that $\mathcal{M}$ edges are in $D_3(\mathcal{F},\mathcal{M})$. As $D_1(\mathcal{F},\mathcal{M})$ edges are counted twice and $D_2(\mathcal{F},\mathcal{M})$ edges are counted once in $\stackrel [{x \in (X_{\mathcal{F}}\setminus X_{\mathcal{M}})}]{}{\sum}|\mathcal{F}_x|$, we get $\stackrel [{x \in (X_{\mathcal{F}}\setminus X_{\mathcal{M}})}]{}{\sum}|\mathcal{F}_x|= 2d_1+d_2$. Similarly, $\stackrel[{x \in X_{\mathcal{M}}}]{}{\sum}|\mathcal{F}_x|=d_1+2d_2+3d_3$. By case assumption and two immediate previous statements, $2d_1+d_2>d_1+2d_2+3d_3$. Therefore, $2d_1-2d_3>d_1+d_2+d_3=|\mathcal{F}|$ as
$\{D_i(\mathcal{F},\mathcal{M})|\ i \in \{1,2,3\}\}$ is a partition of $\mathcal{F}$. Thus,
\begin{eqnarray*}
|\mathcal{F}|&<& 2d_1-2d_3 \\
                         &\leq& 2d_1-2\nu \mbox{ [as $d_3\geq \nu$]}\\
                         &\leq& 2\max\{(d-1)\nu, 6\nu\}-2\nu \mbox{ [as by Lemma \ref{D1 edges} $d_1\leq\max\{(d-1)\nu, 6\nu\}$]}\\
                         &=& 2\nu \max\{(d-2),5\}. 
\end{eqnarray*} 
Therefore,                            
\begin{equation}\label{case_2_eq}
 2\nu \max\{(d-2),5\}\geq |\mathcal{F}|.
\end{equation}
By equations (\ref{case_1_eq}) and (\ref{case_2_eq}), $|\mathcal{F}|\leq \max\{2d\nu, 10\nu\}$. \hspace{\stretch{1}} $\square$\newline

It is challenging to prove our main result--Theorem \ref{thm1}. In the next section some tools are built to prove Theorem \ref{thm1}.
\section{Important Propositions}
To state these useful propositions precisely, we need more notions such as the set of vertices that are covered by each maximum matching. 
\begin{defn}\label{SF} Let $\mathcal{F}$ be a set system. Then $S_{\mathcal{F}}$ denotes the set of vertices, in $X_{\mathcal{F}}=\stackrel[A\in \mathcal{F}]{}{\cup}A$, that are covered by each maximum matching.\end{defn}
Removal of vertices in $S_{\mathcal{F}}$ along with edges containing these vertices has been a crucial step in finding the bound on the cardinality of an edge set of simple graphs in \cite{bk}. We shall use similar ideas in the proceeding work. The following lemma, which is an easy consequence of Lemma \ref{augmenting set lem}, is left for readers to prove.
\begin{lem}\label{matching size SF}
Let $\mathcal{F}$ be a set system and $x \in X_{\mathcal{F}}$. $x\in S_{\mathcal{F}}$ if and only if\\ $\nu(\mathcal{F}\setminus \mathcal{F}_x)=\nu(\mathcal{F})-1$.\end{lem}
We make the following crucial remark based on the lemma above. This remark is one of the key ideas that prove the main result.
\begin{rem}\label{SF-rem} Let $\mathcal{F}$ be a set system with $x \in S_{\mathcal{F}}$. Then
$|\mathcal{F}|\leq |\mathcal{F}_x|+|\mathcal{F}\setminus \mathcal{F}_x|\leq  \Delta(\mathcal{F})+|\mathcal{F}\setminus \mathcal{F}_x|$ and  by Lemma \ref{matching size SF}, $\nu(\mathcal{F}\setminus \mathcal{F}_x)=\nu(\mathcal{F})-1$.
\end{rem}
\begin{defn}\label{nested-sequence}
Let $\mathcal{H}$ be a $k$-uniform family with $S_{\mathcal{H}}\neq \emptyset$. A sequence $(x_1,x_2,\ldots,x_{k_1})$ of verticies of $\mathcal{H}$ is called nested if there exists a corresponding sequence of subfamilies $\mathcal{H}_0$, $\mathcal{H}_1$, \ldots, $\mathcal{H}_{k_1}$ such that $x_i$'s and $\mathcal{H}_i$'s satisfy:\newline
(i) $\mathcal{H}_0:=\mathcal{H}$,\newline
(ii) $x_i\in S_{\mathcal{H}_{i-1}}$ and $\mathcal{H}_i:=\mathcal{H}_{i-1}\setminus \mathcal{H}_{x_i}$ for all $1\leq i\leq k_1$. The positive integer $k_1$ is such that $S_{\mathcal{H}_{k_1}}=\emptyset$.
\end{defn}
Note that the value of $k_1$, defined by \ref{nested-sequence}, depends on the sequence $(x_i)$ for $i \in \{1,\ldots ,k_1\}$ as shown in the example below.
\begin{rem} Let $G$ be the following graph. $V(G)=\{w,x,y,z\}$ and $E(G)=\{\{w,x\},\{x,y\},\{y,z\},\{x,z\}\}$. Note $\{\{w,x\},\{y,z\}\}$ is the only maximum matching of $G$ and hence every vertex is covered by all maximum matchings of $G$. Thus, $S_G=V(G)$ by Definition \ref{SF}. Consider two sequences of vertices $(w)$ and $(x, y)$ for $x_i$'s in the definition \ref{nested-sequence}.
\begin{itemize}
\item[(i)] Let $x_1=w$ and consider induced subgraph $G_1$ on $V(G)\setminus \{w\}$. Then $E(G_1)=E(G)\setminus \{\{w,x\}\}$. Note that any of the three edges of $G_1$, $\{\{x,y\},\{y,z\},\{x,z\}\}$, is a maximum matching of $G_1$. Hence for each vertex $v$ of $G_1$ there is a corresponding maximum matching of $G_1$ not covering $v$ and so $\mathcal{S}_{G_1}=\emptyset$ and $k_1=1$.
\item[(ii)] Let $x_1=x$ and consider induced subgraph $G_2$ on vertices $V(G)\setminus \{x\}$. Then $E(G_2)=E(G)\setminus \{\{w,x\},\{x,y\},\{x,z\}\}=\{y,z\}$. The edge $\{y,z\}$ is the only maximum matching of $G_2$ hence $\{y,z\}\subseteq \mathcal{S}_{G_2}$. In this case $k_1=2$ and any of $y$ or $z$ can be chosen as $x_2$.     
\end{itemize}
\end{rem} 
There are other interesting facts about nested sequences such as reordering of vertices of a nested sequence results in another nested sequence. However, we will not be needing these facts for the following discussion. The lemma below provides a bound on the maximum degree of a $k$-uniform, linear family $\mathcal{F}$ if $S_{\mathcal{F}}=\emptyset$.     
\begin{prop}\label{degree_knu}
Let $\mathcal{F}$ be a $k$-uniform, linear family and let $\nu:=\nu(\mathcal{F})$. If there exists a $x \in X_{\mathcal{F}}$ such that $|\mathcal{F}_x| >k\nu$, then $x \in S_{\mathcal{F}}$.
\end{prop} 
\noindent{\sl Proof.} By Definition \ref{SF}, a vertex $x\in S_{\mathcal{F}}$ if and only if $x$ is covered by every maximum matching of $\mathcal{F}$. Assume on the contrary that $x \notin S_{\mathcal{F}}$. Then there exists a maximum matching $\mathcal{M}$ of $\mathcal{F}$ such that $x \notin X_{\mathcal{M}}$. For any $A \in \mathcal{F}_x$, $A\cap X_{\mathcal{M}}\neq \emptyset$ as $\mathcal{M}$ is a maximum matching. Otherwise there is an $\mathcal{M}$-augmenting set $\{A\}$. However, $\mathcal{F}_x$ is a linear family such that $\stackrel[A \in \mathcal{F}_x]{}{\cap}A=\{x\}$. Thus for any $\{A,B\}\subseteq \mathcal{F}_x$, $(A \cap X_{\mathcal{M}})\cap (B \cap X_{\mathcal{M}})=\emptyset$. Hence $k\nu=|X_{\mathcal{M}}|\geq |X_{\mathcal{F}_x} \cap X_{\mathcal{M}}|\geq|\mathcal{F}_x|$ but this contradicts $|\mathcal{F}_x| >k\nu=|X_{\mathcal{M}}|$.\hspace{\stretch{1}} $\square$\newline

\begin{prop}\label{sf-bound} Let $\mathcal{F}_i$, $x_i$ and $k_1$ be defined as in Definition \ref{nested-sequence}. If $d=\Delta(\mathcal{F})$, then
\begin{equation}\label{eq-sf-bound}
 |\mathcal{F}|\leq k_1d+|\mathcal{F}_{k_1}|.
\end{equation}
Furthermore if $\mathcal{F}$ is a $k$-uniform, linear family then $\Delta({\mathcal{F}}_{k_1})\leq \min\{k\nu({\mathcal{F}}_{k_1}),d\}$.
\end{prop} 
\noindent{\sl Proof.} The equation (\ref{eq-sf-bound}) obviously holds as $|\mathcal{F}_{x_i}|\leq \Delta(\mathcal{F})=d$ for each $i \in \{1,\ldots,k_1\}$ and $\mathcal{F}=\cup^{k_1}_{i=1} (\mathcal{F}_{x_i})\cup\mathcal{F}_{k_1}$. By Proposition \ref{degree_knu}, $\Delta({\mathcal{F}}_{k_1})\leq k\nu({\mathcal{F}}_{k_1})$ or else $\mathcal{S}_{{\mathcal{F}}_{k_1}}\neq \emptyset$ contrary to the definition of $k_1$. Also, $\Delta({\mathcal{F}}_{k_1})\leq \Delta(\mathcal{F})=d$ as ${\mathcal{F}}_{k_1}\subseteq \mathcal{F}$.\hspace{\stretch{1}} $\square$\newline

We next partition $\mathcal{F}$ to establish some crucial propositions. Let $\mathcal{F}$ be a $3$-uniform, linear family,  $\mathcal{M}$ be a maximum matching of $\mathcal{F}$ with $S_{\mathcal{F}}= \emptyset$, $d:=\Delta(\mathcal{F})$ and $\nu:=\nu(\mathcal{F})$. By Proposition \ref{degree_knu}, $d\leq 3\nu$.
Now define as before,
\begin{defn}\label{cij}
Let $\mathcal{F}$ and $\mathcal{M}$ be as described above.\newline
For $i\in \{1,2,3\}$, $D_i(\mathcal{F})=\{A\in \mathcal{F}\quad|\quad |A\cap X_{\mathcal{M}}|=i\}$.\newline
For $i\in \{1,2,3\}$ and $y \in X_{\mathcal{F}}$, $d_i(y)=|\ D_i(\mathcal{F})\cap \mathcal{F}_y|$.\newline
For $A, B \in \mathcal{M}$, $D_2(A,B)= \{C\in D_2(\mathcal{F})\quad|\quad C \cap A\neq \emptyset, C \cap B \neq \emptyset \}$.\newline
For $A,B,C\in \mathcal{M}$, $D_2(A,B,C)= \{E\in D_2(\mathcal{F})\quad|\quad |E \cap (A\cup B\cup C)|=2 \}$.\newline
\end{defn}

\medskip

Note that $\{D_i(\mathcal{F})\quad|\quad i\in\{1,2,3\}\}$ is a partition of $\mathcal{F}$ and $\mathcal{M}\subseteq D_3(\mathcal{F})$. Next, we find bounds on $|D_2(A,B)|$ and $|D_2(A,B,C)|$.  
\begin{prop}\label{prop_c00}
For all $\{A,B\}\subseteq \mathcal{M}$, $|D_2(A,B)|\leq 8$.
\end{prop}
\noindent{\sl Proof.} Let $D(A,B):=\{C\in \mathcal{F}\quad|\quad C\cap A\neq \emptyset,C\cap B\neq \emptyset\}$. Clearly, $D_2(A,B)\subseteq D(A,B)$. Since $\mathcal{F}$ is linear, there is at most one edge of $\mathcal{F}$ that contains both $a$ and $b$ for any $a\in A$ and $b\in B$. Therefore, $D(A,B)\leq 9$. In particular, $D_2(A,B)\leq 9$. Assume $D_2(A,B)= 9$; we shall obtain a contradiction to the fact that $\mathcal{M}$ is a maximum matching.\newline
Let $A=\{1,2,3\}$ and $B=\{4,5,6\}$. We construct a graph $G$ with vertex set $V(G)=\{1,2,3,4,5,6\}$ and edge set $\{\{i,j\}|\ i \in A, j\in B\}$. Since $D_2(A,B)\subseteq D_2(\mathcal{F})$, the only edges of $\mathcal{M}$ covered by edges in $D_2(A,B)$ are $A$ and $B$. Hence if $\{i,j,u\}\in D_2(\mathcal{F})$ with $\{i,j\} \in E(G)$ then $u\notin X_{\mathcal{M}}$. Now consider any matching $N$ of size three in $G$. Without loss of generality, let $N= \{\{1,4\},\{2,5\},\{3,6\}\}$ and let the edges in $D_2(A,B)$ covering $N$ be $\{1,4,u\},\{2,5,v\},\{3,6,w\}$. If no two of $u$, $v$ and $w$ are the same vertex then we have an augmenting set $\{\{1,4,u\},\{2,5,v\},\{3,6,w\},A,B\}$  in $\mathcal{F}$ and $\mathcal{M}$ is not a maximum matching by Lemma \ref{augmenting set lem}. So without loss of generality, let $ v=w$.\newline
{\bf{Claim}}: Let $\{1,4,u\},\{2,5,v\},\{3,6,v\},\{2,4,s\}$, $\{1,5,t\}$, $\{1,6,y\}$ and $\{3,4,z\}$ be edges in $\mathcal{F}$. Then $s=t$ and $y=z$.\newline
{\bf{Proof of the claim}}: Note that $s\neq v$ and $t\neq v$ as the sets $\{2,v\}$ and $\{5,v\}$ are contained in a unique element of $\mathcal{F}$. So, if $s\neq t$ then \newline $\{\{2,4,s\}, \{1,5,t\}, \{3,6,v\}, A, B\}$ is an $\mathcal{M}$-\textit{augmenting set}. But this is a contradiction as $\mathcal{M}$ is a maximum matching and Lemma \ref{augmenting set lem} implies that $\mathcal{F}$ has no $\mathcal{M}$-augmenting set. Also, $y\neq v$ and $z \neq v$ because $\{3,v\}$ and $\{6,v\}$ are contained in a unique element of $\mathcal{F}$. So, if $y\neq z$ then \newline $\{\{1,6,y\}, \{3,4,z\}, \{2,5,v\}, A, B\}$ is an $\mathcal{M}$-\textit{augmenting set} again leading to a contradiction by Lemma \ref{augmenting set lem}. Thus, the claim is established.\newline
If $\{2,6,r\}\in \mathcal{F}$ then $r\neq y$ because $\{1,6,y\}\in \mathcal{F}$ contains $\{6,y\}$ and $r\neq s$ because $\{2,4,s\}\in \mathcal{F}$ contains $\{2,s\}$. Hence the above claim implies that \newline
$\{\{1,5,s\},\{3,4,y\},\{2,6,r\},A,B\}$ is an $\mathcal{M}$-\textit{augmenting set}, leading to a contradiction by Lemma \ref{augmenting set lem}.\hspace{\stretch{1}} $\square$\newline
\begin{rem}\label{rem_c00}
Up to isomorphism, there exists a unique configuration of eight edges in $D_2(A,B)$. Namely, if  $A=\{1,2,3\}$ and $B=\{4,5,6\}$ then $D_2(A,B)$ is: $\{\{1,5,s\},\{2,6,s\},\{1,4,t\},\{3,5,t\},\{2,4,u\},\{1,6,u\},\{2,5,v\},\{3,4,v\}\}$\newline where $s$, $t$, $u$ and $v$ are different vertices. Readers can establish the uniqueness by arguing similarly to Proposition \ref{prop_c00}.
\end{rem}
Now we find the maximum value of $|D_2(A,B,C)|$. It is clear that\newline $|D_2(A,B,C)|\leq |D_2(A,B)|+|D_2(B,C)|+|D_2(A,C)|\leq 24 $. We improve the bound to $|D_2(A,B,C)|\leq 21$ in the next two propositions. 
\begin{defn}\label{g2}
Let $\mathcal{F}$ be a $3$-uniform, linear family, and let $\mathcal{M}$ be a matching (need not be maximum) of $\mathcal{F}$. For  $\{A,B\}\subseteq \mathcal{M}$, we define a simple graph $G(D_2,A,B)$ as follows:
$V(G(D_2,A,B)):= A\cup B$ and $E(G(D_2,A,B)):=\{C \cap (A \cup B)\quad|\quad C\in D_2(A,B)\}$.
\end{defn}
\begin{prop}\label{prop1 c000} Let $\mathcal{F}$ be a linear, $3$-uniform family and let $\mathcal{M}$ be a maximum matching of $\mathcal{F}$. If $\{A, B, C\}\subseteq \mathcal{M}$ and $|D_2(A,B)|=8$ then $|D_2(A,C)|+|D_2(B,C)|\leq 12$.
\end{prop}
\noindent{\sl Proof.} Let $A=\{1,2,3\}$, $B=\{4,5,6\}$ and $C=\{7,8,9\}$. As $D_2(A,B)=8$, without loss of generality let $\{3,6\}\notin E(G(D_2,A,B))$ and hence\newline $E(G(D_2,A,B))=\{\{1,4,\},\{1,5\},\{1,6\},\{2,4\},\{2,5\},\{2,6\},\{3,4\},\{3,5\}\}$. Also without loss of generality, by Remark \ref{rem_c00}, the subfamily corresponding to $G(D_2,A,B)$ in $\mathcal{F}$ is \begin{equation}{\{\{1,5,s\},\{2,6,s\},\{1,4,t\},\{3,5,t\},\{2,4,u\},\{1,6,u\},\{2,5,v\},\{3,4,v\}\}\label{unique_family_c00}}\end{equation} where $s$, $t$, $u$ and $v$ are different vertices and are not covered by the maximum matching $\mathcal{M}$.\newline
\begin{claim}\label{claim1_propc000}: $|\{E\in D_2(\mathcal{F})\quad |\quad |E\cap \{7,8,9\}|=1,\ |E\cap \{3,6\}|=1\}|\leq 4$.\end{claim}
{\bf{Proof of Claim \ref{claim1_propc000}}}: If the claim does not hold then without loss of generality $3$ edges of $D_2(A,C)$ are incident to the vertex $3$ and at least $2$ edges of $D_2(B,C)$ are incident to the vertex $6 $. We may assume that there are edges $\{6,7,a\}$ and $\{6,8,b\}$ in $D_2(B,C)$. 
By our assumption (\ref{unique_family_c00}),\newline
$\{\{1,5,s\},\{2,6,s\},\{1,4,t\},\{3,5,t\},\{2,4,u\},\{1,6,u\},\{2,5,v\},\{3,4,v\}\}\subseteq \mathcal{F}$ where $s$, $t$, $u$ and $v$ are different vertices and are not covered by the maximum matching $\mathcal{M}$. Also by assumption $\{\{3,7,x\},\{3,8,y\},\{3,9,z\},\{6,7,a\},\{6,8,b\}\}\subseteq \mathcal{F}$ for some vertices $x$, $y$, $z$, $a$ and $b$ in $X_{\mathcal{F}}\setminus X_{\mathcal{M}}$.\newline
We will use the following two observations.\newline
(i) As $\{\{3,7,x\},\{3,8,y\},\{3,9,z\}, \{3,4,v\},\{3,5,t\}\}\subseteq \mathcal{F}$ and $\mathcal{F}$ is a linear family, $x \notin \{t,v\}$, $y\notin \{t,v\}$ and $z\notin \{t,v\}$.\newline
(ii) As $\{ \{2,6,s\},\{1,6,u\}, \{6,7,a\},\{6,8,b\}\}\subseteq \mathcal{F}$ and $\mathcal{F}$ is a linear family, $a\notin \{s,u\}$ and $b\notin\{s,u\}$.\newline
Suppose that $x=b$. Then $z\neq b$ because $\{\{3,7,x\},\{3,9,z\}\}\subseteq \mathcal{F}$. Since $z\neq b$, $z \notin \{t,v\}$ and $b=x\notin \{t,v\}$, we have an $\mathcal{M}$-augmenting set \newline
$\{\{1,4,t\},\{2,5,v\},\{3,9,z\},\{6,8,b\},A,B,C\}$ in $\mathcal{F}$ contradicting Lemma \ref{augmenting set lem} as $\mathcal{M}$ is a maximum matching of $\mathcal{F}$. Symmetrically, if $z=b$ then $x\neq b$ because $\{\{3,7,x\},\{3,,9,z\}\}\subseteq \mathcal{F}$. Since $x \neq b$, $x \notin \{t,v\}$ and $b=z\notin \{t,v\}$, we have an $\mathcal{M}$-augmenting set\newline  $\{\{1,4,t\},\{2,5,v\},\{3,7,x\},\{6,8,b\},A,B,C\}$ in $\mathcal{F}$.\newline
So far we have shown that $b\notin\{x,z\}$. We claim that $\{x,z\}=\{s,u\}$. If this claim doesn't hold then either $x \notin \{s,u\}$ or  $z \notin \{s,u\}$. Let $x \notin \{s,u\}$. The case  $z \notin \{s,u\}$ is similar. Since $x\notin \{s,u\}$, $x\neq b$ and by observation (ii) $b \notin \{s,u\}$, we get the following $\mathcal{M}$-augmenting set\newline $\{\{3,7,x\},\{2,4,u\},\{1,5,s\},\{6,8,b\},A,B,C\}$, a contradiction.\newline
Finally, note that $a \neq z$ as $z \in \{s,u\}$ and by observation (ii) $a\notin \{s,u\}$. Next we claim that $a \in \{t,v\}$. If this claim doesn't hold then\newline $\{\{6,7,a\},\{1,4,t\},\{2,5,v\},\{3,9,z\},A,B,C\}$ is an $\mathcal{M}$-augmenting set. Thus, $a \in \{t,v\}$ and by observation (i) $y\notin \{t,v\}$. Therefore, $a \neq y$. Note that $y \notin \{s,u\}$ as $\{x,z\}=\{s,u\}$. So, we have the following $\mathcal{M}$-augmenting set\newline $\{\{2,4,u\},\{1,5,s\},\{3,8,y\},\{6,7,a\},A,B,C\}$ in $\mathcal{F}$. This contradiction to the maximality of $\mathcal{M}$ completes the proof of Claim \ref{claim1_propc000}.\newline
       
For $i \in \{7,8,9\}$, define $D_2(i):=\{E\in D_2(\mathcal{F})\quad|\quad E\cap \{1,2,4,5\}\neq \emptyset \mbox{ and $i \in E$}\}.$\newline
\begin{claim}\label{claim2_propc000}: For $\{i,j\}\subseteq \{7,8,9\}$, $|D_2(i)|+|D_2(j)|\leq 6$.\end{claim}
{\bf{Proof of Claim \ref{claim2_propc000}}}: Without loss of generality, let $i=7$ and $j=8$ and assume on the contrary $|D_2(7)|+|D_2(8)|\geq 7$. As $|D_2(i)|\leq 4$ for $i\in \{7,8,9\}$ by definition, without loss of generality let $|D_2(7)|=4$ and $|D_2(8)|\geq 3$. Also by symmetry of $1$, $2$, $4$, $5$, we may assume that there are edges in $D_2(7)\cup D_2(8)$ containing each of $\{\{1,7\},\{1,8\},\{2,7\},\{2,8\},\{4,7\},\{4,8\},\{5,7\}\}$.  By our initial assumption (\ref{unique_family_c00}),\newline
$\{\{1,5,s\},\{2,6,s\},\{1,4,t\},\{3,5,t\},\{2,4,u\},\{1,6,u\},\{2,5,v\},\{3,4,v\}\}\subseteq \mathcal{F}$ where $s$, $t$, $u$ and $v$ are different vertices and are not covered by the maximum matching $\mathcal{M}$. Let\newline $\{\{1,7,a\},\{2,7,b\},\{4,7,c\},\{5,7,d\},\{1,8,x\},\{2,8,y\},\{4,8,z\}\}\subseteq \mathcal{F}$. The $\{0,1\}$-intersection property of $\mathcal{F}$ implies that $a \notin \{t,s,u,b,c,d\}$, $b \notin \{s,v,u,a,c,d\}$, $c\notin \{t,v,u,a,b,d\}$, $d\notin \{t,s,v,a,b,c\}$, $x\notin \{s,t,u,y,z,a\}$, $y\notin\{s,v,u,x,z,b\}$ and $z \notin \{t,u,v,x,y,c\}$. We now make observations that prove Claim \ref{claim2_propc000}.\newline
\begin{fact}\label{fact_2_1}
Either $c=x$ or $c=s$.
\end{fact}
{\sl Proof.} We have $t\neq s$, $c\neq t$ and $x \notin \{t,s\}$. If $c\notin \{x,s\}$ then\newline
 $\{\{4,7,c\},\{1,8,x\},\{3,5,t\},\{2,6,s\},A,B,C\}$ is an $\mathcal{M}$-augmenting set in $\mathcal{F}$, a contradiction.
\begin{fact}\label{fact_2_2}
$b=t$.
\end{fact}
{\sl Proof.} Since $t\neq u$, $z \notin \{t,u\}$ and $b\neq u$, either $b=t$ or $b=z$ otherwise $\{\{2,7,b\},\{3,5,t\},\{4,8,z\},\{1,6,u\},A,B,C\}$ is an $\mathcal{M}$-augmenting set in $\mathcal{F}$.\newline
If $b=z$ then $b \notin \{s,t,u,v,x,y,a,c,d\}$ as noted earlier. But then we have the following $\mathcal{M}$-augmenting set\newline
$\{\{1,7,a\},\{2,6,s\},\{3,5,t\},\{4,8,b\},A,B,C\}$ in $\mathcal{F}$, a contradiction.
\begin{fact}\label{fact_2_3}
$y=t$.
\end{fact}
{\sl Proof.} Since $t\neq u$, $c\notin \{t,u\}$ and $y\neq u$, either $y=t$ or $y=c$ otherwise $\{\{2,8,y\},\{1,6,u\},\{4,7,c\},\{3,5,t\},A,B,C\}$ is an $\mathcal{M}$-augmenting set in $\mathcal{F}$. But $c\neq y$ because $c \in \{s,x\}$ by Fact \ref{fact_2_1} and, as noted prior to Fact \ref{fact_2_1}, $y \notin \{s,x\}$. This completes the proof of this Fact.\newline
\medskip
By Fact \ref{fact_2_2} and Fact \ref{fact_2_3}, $y=t=b$. But this contradicts linearity of the family $\mathcal{F}$ as $|\{2,7,t\}\cap\{2,8,t\}|=2$ and proves Claim \ref{claim2_propc000}. \newline

\medskip

The above claim implies that there can't be strictly more than nine $D_2(\mathcal{F})$ edges such that each edge covers a vertex in $\{1,2,4,5\}$ and another in $\{7,8,9\}$. The next claim improves the estimate. Note that by Claim \ref{claim2_propc000}, if $D_2(i)=4$ for any $i \in \{7,8,9\}$ then $D_2(j)\leq 2$ for $j \in \{7,8,9\}\setminus \{i\}$. Note also that $D_2(i)\leq 4$ by definition and linearity of $\mathcal{F}$.\newline
\begin{claim}\label{claim3_propc000}: There can't be nine or more $D_2(\mathcal{F})$ edges such that each edge covers a vertex in $\{1,2,4,5\}$ and another in $\{7,8,9\}$.\end{claim}
{\bf{Proof of Claim \ref{claim3_propc000}}}: We already know by Claim \ref{claim2_propc000} that there can't be strictly more than nine edges satisfying the condition in Claim \ref{claim3_propc000}. If there are nine such edges then each vertex in $\{7,8,9\}$ is incident to exactly three of $\{1,2,4,5\}$ or Claim \ref{claim2_propc000} is contradicted.\newline
We consider the bipartite graph $G$ on vertices $\{\{1,2,4,5\},\{7,8,9\}\}$ defined by edges in $D_2(A,C)\cup D_2(B,C)$. For all $i\in \{7,8,9\}$, we have $d_G(i)=3$. Since $\lceil{\frac{9}{4}}\rceil=3$, there is a vertex of degree at least three in $\{1,2,4,5\}$. Without loss of generality, we may assume that $d_G(1)\geq3$; the cardinality of the class $\{7,8,9\}$ imposes that $d_G(1)=3$ and that the vertex $1$ is a neighbor of each vertex in $\{7,8,9\}$. Since $d_G(4)+d_G(5)\geq 9-d_G(1)-d_G(2)\geq 3$, either $d_G(4)\geq 2$ or $d_G(5)\geq 2$. So, without loss of generality, let $d_G(4)\geq 2$. Also we can assume that $\{4,7\}$ and $\{4,8\}$ are in $E(G)$ (if not, then reorder vertices $7$, $8$ and $9$). Hence $\{\{1,7,a\},\{1,8,b\},\{1,9,c\},\{4,7,x\},\{4,8,y\}\}\subseteq \mathcal{F}$ for some $a$, $b$, $c$, $x$ and $y$ in $X_{\mathcal{F}}\setminus X_{\mathcal{M}}$. And by our assumption (\ref{unique_family_c00}),\newline $\{\{1,5,s\},\{2,6,s\},\{1,4,t\},\{3,5,t\},\{2,4,u\},\{1,6,u\},\{2,5,v\},\{3,4,v\}\}\subseteq \mathcal{F}$ where $s$, $t$, $u$ and $v$ are different vertices and are not covered by the maximum matching $\mathcal{M}$. The $\{0,1\}$-intersection property implies that $a\notin\{b,c,x,s,t,u\}$, $b\notin \{a,c,y,s,t,u\}$, $c\notin \{a,b,s,t,u\}$, $x\notin \{y,a,t,u,v\}$ and $y\notin \{x,b,t,u,v\}$.
\begin{fact}\label{fact_12_3_1} $x=s$.
\end{fact}
{\sl Proof.} We have $t\neq s$, $c\notin \{t,s\}$ and $x\neq t$. If $x\notin\{c,s\}$, then\newline
$\{\{1,9,c\},\{4,7,x\},\{3,5,t\},\{2,6,s\},A,B,C\}$ is an $\mathcal{M}$-augmenting set in $\mathcal{F}$, a contradiction.
If $x=c$, then $c=x \notin\{a,b,s,t,u,v,y\}$ as noted before Fact \ref{fact_12_3_1}. We also know that $b \notin \{s,t\}$. But then we have the following $\mathcal{M}$-augmenting set $\{\{1,8,b\},\{4,7,x\},\{3,5,t\},\{2,6,s\},A,B,C\}$ in $\mathcal{F}$, a contradiction. Hence $x=s$.
\begin{fact}\label{fact_12_3_2} $y=s$.
\end{fact}
{\sl Proof.} We have $t\neq s$, $c\notin \{t,s\}$ and $y\neq t$. If $y\notin\{c,s\}$, then\newline
$\{\{1,9,c\},\{4,8,y\},\{3,5,t\},\{2,6,s\},A,B,C\}$ is an $\mathcal{M}$-augmenting set in $\mathcal{F}$, a contradiction.
If $y=c$, then $c=y\notin\{a,b,x,s,t,u,v\}$ as noted prior to the previous fact. But this gives the following $\mathcal{M}$-augmenting set\newline $\{\{1,7,a\},\{4,8,y\},\{3,5,t\},\{2,6,s\},A,B,C\}$ in $\mathcal{F}$. Thus, contradicts that $\mathcal{M}$ is a maximum matching. 

By Fact \ref{fact_12_3_1} and Fact \ref{fact_12_3_2}, $x=y=s$. But this contradicts the linearity of $\mathcal{F}$ as $|\{4,7,s\}\cap\{4,8,s\}|=2$. Hence, Claim \ref{claim3_propc000} is proved.
      
The statement of Proposition \ref{prop1 c000} is an easy consequence of Claim \ref{claim1_propc000} and Claim \ref{claim3_propc000}.\hspace{\stretch{1}} $\square$\newline

We shall not be using the following remark. Though, the statement of the remark can improve the bound in the main result as done in author's doctoral dissertation \cite{K}. However, the statement below was proved using the aid of a computer program and we decided not to use it for the current article since the improvement in the bound is not significant. Using the remark below, it can be shown that $|D_2(A,B,C)|\leq 20$ in Proposition \ref{prop_c000}.   
\begin{rem}\label{rem_777} Let $\mathcal{F}$ be a $3$-uniform, linear family and $\mathcal{M}$ be a maximum matching of $\mathcal{F}$. If $\{A,B,C\}\subseteq \mathcal{M}$, then $|D_2(A,B)|=|D_2(A,C)|=D_2(B,C)|=7$ doesn't hold. 
\end{rem}
\begin{prop}\label{prop_c000}
Let $\mathcal{F}$ be a $3$-uniform, linear family and let $\mathcal{M}$ be a maximum matching of $\mathcal{F}$. If $\{A,B,C\}\subseteq \mathcal{M}$, then $|D_2(A,B,C)|\leq 21$. 
\end{prop}
{\sl Proof.} Assume on the contrary \resizebox{.6\hsize}{!}{$|D_2(A,B)|+|D_2(B,C)|+|D_2(A,C)|=|D_2(A,B,C)|\geq 22$}. Therefore, by Proposition \ref{prop_c00} at least one of $|D_2(A,B)|$, $|D_2(B,C)|$ or $|D_2(A,C)|$ is equal to $8$. Without loss of generality, let $D_2(A,B)=8$. Thus, $|D_2(B,C)|+|D_2(A,C)|\geq 13$. This contradicts Proposition \ref{prop1 c000}.\hspace{\stretch{1}} $\square$
\section{$3$-uniform, linear families $\mathcal{F}$ with $S_{\mathcal{F}}= \emptyset$}

In this section, we find a bound on the size of $3$-uniform, linear families $\mathcal{F}$ with $S_{\mathcal{F}}=\emptyset$ (defined by \ref{SF}) in terms of their maximum matching and maximum degree. The chief idea of the proof that establishes the bound follows. For a $3$-uniform, linear family with $\Delta$ approximately greater than  $4\nu$, if $|\mathcal{F}|>\Delta \nu$ then for any given maximum matching $\mathcal{M}$, a local augmenting set involving at most three matching edges is found and extended to a global $\mathcal{M}$-augmenting set. Thus, contradicting the fact that $\mathcal{M}$ is a maximum matching and so establishing the result.
   
Let us recall few notations. Let $\mathcal{F}$ be a $3$-uniform, linear family, and let $\mathcal{M}$ be a maximum matching of $\mathcal{F}$. For $A\in \mathcal{M}$, define
$D_1(A):=\{B \in D_1(\mathcal{F}, \mathcal{M})\quad|\quad B \cap A \neq \emptyset\}$ and $d_1(A):=|D_1(A)|$.
For any $\mathcal{G}\subseteq \mathcal{F}$ and $A \in \mathcal{F}$, also define $\mathcal{G}_A:=\{B\in \mathcal{G}\quad|\quad B\cap A \neq \emptyset\}$.

The following partition of a maximum matching is crucial to obtain the bound on the size of a $3$-uniform, linear family. 
\begin{defn}\label{rem matching partition}
Let $\mathcal{F}$ be a $3$-uniform, linear family with $S_{\mathcal{F}}=\emptyset$, $\nu :=\nu(\mathcal{F})$, $\Delta:=\Delta(\mathcal{F})$ and let $\mathcal{M}$ be a maximum matching of $\mathcal{F}$. We partition $\mathcal{M}$ the following way.\newline
$\mathcal{M}_1:=\{A\in \mathcal{M}\quad|\quad d_1(A)\geq 7\}$ and 
$\mathcal{M}_2:=\mathcal{M}\setminus \mathcal{M}_1$. Also let $m:=|\ \mathcal{M}_1|$ and $\mathcal{M}_1=\{A_1,\ldots,A_m\}$. We already know by Lemma \ref{degree 0-1} that if for some $A \in \mathcal{M}$, $d_1(A)\geq 7$ then all edges in $ D_1(A)$ are incident to the same vertex of $A$. For each $i \in \{1,\ldots,m\}$, let this unique vertex be denoted by $x_i\in A_i$ and let $A_i=\{x_i,y_i,z_i\}$. 
\end{defn}
Since $S_{\mathcal{F}}=\emptyset$, Proposition \ref{degree_knu} implies that $\Delta \leq 3\nu$. Let $\mathcal{M}$, $\mathcal{M}_1$, $\mathcal{M}_2$, $A_i$'s, $x_i$'s, $y_i$'s and $z_i$'s be as defined in the previous definition. Let us partition the family $\mathcal{F}$ and obtain bounds on the size of each class. Since an arbitrary maximum matching $\mathcal{M}$ is fixed in the following discussion, for all $i\in \{1,2,3\}$ $D_i(\mathcal{F})$ is used instead of $D_i(\mathcal{F},\mathcal{M})$.
\begin{defn}\label{partition of family} Let the family $\mathcal{F}$ and $\mathcal{M}$ be as stated in Definition \ref{rem matching partition}. We define\newline
$\mathcal{E}_1:=\stackrel[{i \in \{1,\ldots,m\}}]{}{\cup}\mathcal{F}_{x_i}$;\newline
$\mathcal{E}_2:=\{A \in \mathcal{F}\quad |\quad (\mathcal{M}_2)_A=\emptyset\}\setminus \mathcal{E}_1$, where \newline
$(\mathcal{M}_2)_A= \{B \in \mathcal{M}_2\quad|\quad A \cap B \neq \emptyset\}$, i.e., $\mathcal{E}_2$ consists of those $D_2(\mathcal{F})$ and $D_3(\mathcal{F})$ edges which do not intersect matching edges from $\mathcal{M}_2$ and do not contain vertices from $\{x_1,\ldots,x_m\}$. Note that if $B \in D_1(\mathcal{F})$ then $B\cap (\{y_1,\ldots, y_m\}\cup \{z_1,\ldots,z_m\})= \emptyset$ by Definition \ref {rem matching partition}; \newline
$\mathcal{E}_3:=\{A\in \mathcal{F}\quad|\quad |A \cap X_{\mathcal{M}_2}|=1\}\setminus \mathcal{E}_1$;\newline
$\mathcal{E}_4:=(\{A\in \mathcal{F}\quad |\quad |A \cap X_{\mathcal{M}_2}|\geq 2\}\setminus \mathcal{E}_1)\setminus \mathcal{M}_2$.
\end{defn} 
\begin{rem}\label{rem partition of family}
By Definition \ref{partition of family}, it is obvious that $\mathcal{F}=\cup_{i\in\{1,\ldots,4\}}\mathcal{E}_i \cup \mathcal{M}_2$ and the sets are pairwise disjoint.
\end{rem}
Next we find an upper bound for each member in the above partition with $m=|\mathcal{M}_1|$.
\begin{prop}\label{prop bf-1}
If $\mathcal{E}_1$ is defined by Definition \ref{partition of family}, then $|\mathcal{E}_1|\leq m\Delta$.
\end{prop}
\noindent{\sl Proof.} This is obvious as $\mathcal{E}_1=\cup_{i \in \{1,\ldots,m\}}\mathcal{F}_{x_i}$ and $|\mathcal{F}_{x_i}|\leq \Delta$ for all $i\in \{1,\ldots,m\}$.\hspace{\stretch{1}} $\square$
\begin{prop}\label{prop bf-2}
If $\mathcal{E}_2$ is defined by Definition \ref{partition of family}, then $|\mathcal{E}_2|=0$.
\end{prop}
\noindent{\sl Proof.} Suppose $\mathcal{E}_2\neq \emptyset$, then there exists an edge $B\in \mathcal{E}_2$. By the note after the definition of $\mathcal{E}_2$ (Definition \ref{partition of family}), $B \in D_2(\mathcal{F})\cup D_3(\mathcal{F})$ and all vertices in $B\cap X_{\mathcal{M}}$ belong to $\{y_1,\ldots,y_k\}\cup \{z_1,\ldots,z_k\}$. We show that if $B\in D_2(\mathcal{F})$ or $B \in D_3(\mathcal{F})$, then an $\mathcal{M}$-augmenting set exists in $\mathcal{F}$. Suppose $B \in D_2(\mathcal{F})$. Without loss of generality, let $\{y_1,y_2\}\subseteq B$ and $B=\{y_1,y_2,w\}$ where $w \notin X_{\mathcal{M}}$. Since at least seven $D_1(\mathcal{F})$ edges are incident to $x_1$, at least other seven $D_1(\mathcal{F})$ edges are incident to $x_2$, and there can be at most one edge containing both $w$ and $x_i$ for each $i\in \{1,2\}$, there is an $\mathcal{M}$-augmenting set which consists of an edge from $D_1(\mathcal{F})\cap \mathcal{F}_{x_1}$, an edge from $D_1(\mathcal{F})\cap \mathcal{F}_{x_2}$, B, $\{x_1,y_1,z_1\}$ and $\{x_2,y_2,z_2\}$. This contradicts that $\mathcal{M}$ is a maximum matching. Also for $B \in D_3(\mathcal{F})\cap \mathcal{E}_2$, we can similarly construct an $\mathcal{M}$-augmenting set in $\mathcal{F}$. In this case the augmenting set consists of three $D_1(\mathcal{F})$ edges, the edge $B$ and the three $\mathcal{M}_1$ edges that have nonempty intersection with $B$. Hence in either case there is an $\mathcal{M}$-augmenting set. Thus, $\mathcal{E}_2=\emptyset$. \hspace{\stretch{1}} $\square$ 
\begin{prop}\label{prop bf-3}
If $\mathcal{E}_3$ is defined by Definition \ref{partition of family}, then\newline $|\mathcal{E}_3|\leq \min \{2m+6,\Delta-1\}(\nu-m)$.
\end{prop}
\noindent{\sl Proof.} Recall that $\mathcal{E}_3=\{A \in \mathcal{M}\quad|\quad |A \cap X_{\mathcal{M}_2}|=1\}\setminus \mathcal{E}_1$. Hence $\mathcal{E}_3$ consists of $D_1(\mathcal{F})$ edges that intersect $\mathcal{M}_2$ edges and $D_2(\mathcal{F})\cup D_3(\mathcal{F})$ edges that cover exactly one vertex in $X_{\mathcal{M}_2}$ and no vertex in $\{x_1,\cdots, x_m\}$.\newline
{\bf{Claim }}: If seven or more edges from $\mathcal{E}_3$ intersect an edge $A\in \mathcal{M}_2$ then all $\mathcal{E}_3$ edges that intersect $A$ must be incident to the same vertex $x$ in $A$.\newline
{\bf{Proof of the claim }}: Suppose not; then there exist $B_1$ and $B_2$ in $\mathcal{E}_3$ that intersect $A$ and are disjoint. As at least seven edges from $\mathcal{E}_3$ intersect $A$ and $|A|=3$, by pigeonhole principal there is a vertex $a \in A$ such that among $\mathcal{E}_3$ edges that intersect $A$ at least three contain $a$. If there exists $B_1 \in \mathcal{E}_3$ such that $ B_1$ intersects $A$ and $a \notin B_1$ then we can choose $B_2$ among the edges in $\mathcal{E}_3$ containing $a$.\newline
If $B_1$ and $B_2$ are both $D_1(\mathcal{F})$ edges then $\{B_1,B_2,A\}$ is an $\mathcal{M}$-augmenting set. Now we consider all remaining possibilities for $B_1$ and $B_2$. Considering symmetries, we have the following possibilities.\newline
(i) $B_1$ is a $D_2(\mathcal{F})$ edge and $B_2$ is a $D_1(\mathcal{F})$ edge;\newline
(ii) $B_1$ is a $D_2(\mathcal{F})$ edge and $B_2$ is a $D_2(\mathcal{F})$ edge;\newline
(iii) $B_1$ is a $D_3(\mathcal{F})$ edge and $B_2$ is a $D_1(\mathcal{F})$ edge;\newline
(iv) $B_1$ is a $D_3(\mathcal{F})$ edge and $B_2$ is a $D_2(\mathcal{F})$ edge;\newline
(v) $B_1$ is a $D_3(\mathcal{F})$ edge and $B_2$ is a $D_3(\mathcal{F})$ edge.\newline
In each of the above cases, an $\mathcal{M}$-augmenting set can be constructed using $D_1(\mathcal{F})$ edges incident at $\mathcal{M}_1$ edges along with the $\mathcal{M}_1$ edges intersected by $B_1$ and $B_2$, $B_1$, $B_2$ and $A$. For example consider the case (v), since $B_1$ and $B_2$ are in $D_3(\mathcal{F})$ each of them covers two edges from $\mathcal{M}_1$. Assume the worst case that $B_1$, $B_2$ intersect four different edges in $\mathcal{M}_1$ and let the edges be $A_1$, $A_2$, $A_3$ and $A_4$. Recall that seven or more $D_1(\mathcal{F})$ edges are incident to $x_i\in A_i$ for all $i\in \{1,\ldots,m\}$. Note that any $D_1(\mathcal{F})$ edge incident at $x_i$ can at most intersect two $D_1(\mathcal{F})$ edges incident at $x_j$ for $i\neq j$. Hence there are four\footnote{ We need at least seven $D_1(\mathcal{F})$ edges to be incident at each of the $x_i$'s to ensure existence of four pairwise disjoint $D_1(\mathcal{F})$ edges.} pairwise disjoint $D_1(\mathcal{F})$ edges in $\cup_{i=1}^4 (D_1(\mathcal{F})\cap \mathcal{F}_{x_i})$. These disjoint edges along with $A_1$, $A_2$, $A_3$, $A_4$, $B_1$, $B_2$ and $A$ form an $\mathcal{M}$-augmenting set, a contradiction.\newline
Hence, if seven or more $\mathcal{E}_3$ edges intersect with any $M_2$ edge then all these edges must contain the same vertex of the $\mathcal{M}_2$ edge. Now consider $({\mathcal{E}_3})_A$, the set of $\mathcal{E}_3$ edges incident at an $\mathcal{M}_2$ edge A. If $|({\mathcal{E}_3})_A\cap D_1(\mathcal{F})|\geq 7$, then all $(D_1(\mathcal{F}))_A$ edges are incident to the same vertex in $A$ and $A\in \mathcal{M}_1$. A contradiction to the fact that $A \in \mathcal{M}_2$. Therefore, there are at most six $D_1(\mathcal{F})$ edges in $({\mathcal{E}_3})_A$. By Definition \ref{partition of family}, an edge in $\mathcal{E}_3$ is either a $D_1(\mathcal{F})$ edge or a $D_2(\mathcal{F})\cup D_3(\mathcal{F})$ edge that contains at least one vertex in $\{y_1,\cdots,y_m\}\cup\{z_1,\cdots,z_m\}$ and no vertex in $\{x_1,\cdots,x_m\}$. Hence $|({\mathcal{E}_3})_A|\leq \min\{2m+6,\Delta-1\}$ for all $A \in \mathcal{M}_2$. Therefore, $|\mathcal{E}_3|\leq \min \{2m+6,\Delta-1\}(\nu-m)$. \hspace{\stretch{1}}$\square$ 

Let us generalize Definition \ref{cij} to find a bound on $D_2(\mathcal{F}, \mathcal{M})\cup D_3(\mathcal{F}, \mathcal{M})$.
\begin{defn}\label{d2_f_m}
Let $\mathcal{F}$ be a $3$-uniform, linear family and let $\mathcal{M}$ be a matching (not necessarily maximum) of $\mathcal{F}$. For $i \in \{0,1,2,3\}$, define\newline
$D_i(\mathcal{F},\mathcal{M}):=\{A \in \mathcal{F}\quad|\quad |A \cap X_{\mathcal{F}}|=i\}$. Also for all $\{A,B,C\}\subseteq  \mathcal{M}$, define
$D_2(A,B,C):= \{E\in D_2(\mathcal{F},\mathcal{M})\quad|\quad |E \cap (A\cup B\cup C)|=2 \}$ and\newline
$D_3(A,B,C):= \{E\in (D_3(\mathcal{F},\mathcal{M})\setminus \{A,B,C\})\quad|\quad |E \cap (A\cup B\cup C)|\geq2 \}$.
\end{defn}
\begin{prop}\label{prop_imp d2+d3}
Let $\mathcal{F}$ be a $3$-uniform, linear family and let $\mathcal{M}$ be a matching (not necessarily maximum) of $\mathcal{F}$ such that $n=|\mathcal{M}|$. If $|D_2(A,B,C)|\leq 21$ for all $\{A,B,C\}\subseteq \mathcal{M}$, then $|D_2(\mathcal{F}, \mathcal{M})|+|D_3(\mathcal{F},\mathcal{M})\setminus \mathcal{M}|\leq \frac{23}{(n-2)}{n \choose 3}$.\newline
\end{prop}
{\sl Proof.} For $\{A,B,C\}\subseteq \mathcal{M}$, let\newline
\resizebox{.99 \hsize}{!}{$\mathcal{H}(A,B,C):=\{\{i,j\}\quad|\quad \{i,j\} \mbox{ is contained in an edge from } D_2(A,B,C)\cup(D_3(A,B,C)\setminus \mathcal{M})\}$}.\newline Since $\mathcal{F}$ is a linear family, we get 
$|\{E \in \mathcal{F}\quad|\quad |E\cap (A \cup B)|= 2\}|\leq 9$ for any $\{A,B\}\subseteq \mathcal{M}$. Thus, we obtain 
\begin{equation}\label{d2+d3_eq_1}
|\mathcal{H}(A,B,C)|\leq 27.
\end{equation}
In the expression 
\begin{equation}\label{eq_star}
\sum_{\{A,B,C\}\subseteq \mathcal{M}}|\mathcal{H}(A,B,C)| 
\end{equation} each edge in $D_2(\mathcal{F},\mathcal{M})$ is counted $(n-2)$ times because $C$ can be any of the $(n-2)$ other $\mathcal{M}$ edges for a fixed pair $\{A,B\}\subset \mathcal{M}$. Also each edge in $D_3(\mathcal{F},\mathcal{M})\setminus \mathcal{M}$ is counted $3(n-2)$ times in the expression (\ref{eq_star}). Hence
\begin{equation}\label{d2+d3_eq_2}
(n-2)|D_2(\mathcal{F},\mathcal{M})|+3(n-2)|D_3(\mathcal{F})\setminus \mathcal{M}|=\sum_{\{A,B,C\}\subseteq \mathcal{M}}|\mathcal{H}(A,B,C)| .\end{equation} So by equations (\ref{d2+d3_eq_1}) and (\ref{d2+d3_eq_2}), we have
\begin{displaymath}
(n-2)|D_2(\mathcal{F},\mathcal{M})|+3(n-2)|D_3(\mathcal{F})\setminus \mathcal{M}|\leq 27{n \choose 3}.
\end{displaymath}
Therefore,
\begin{equation}\label{d2+d3_eq_3}
|D_3(\mathcal{F})\setminus \mathcal{M}|\leq \frac{27}{3(n-2)}{n \choose 3}- \frac{1}{3}|D_2(\mathcal{F},\mathcal{M})|.
\end{equation}
So, we have
\begin{equation}\label{d2+d3_eq_4}
|D_2(\mathcal{F},\mathcal{M})|+|D_3(\mathcal{F})\setminus \mathcal{M}|\leq \frac{2}{3}|D_2(\mathcal{F},\mathcal{M})|+\frac{27}{3(n-2)}{n \choose 3}.
\end{equation}
By equation (\ref{d2+d3_eq_2}), we have
\begin{equation}
(n-2)|D_2(\mathcal{F},\mathcal{M})|=\sum_{\{A,B,C\}\subseteq \mathcal{M}}|\mathcal{H}(A,B,C)\cap D_2(\mathcal{F},\mathcal{M})|.
\end{equation}
As $$\mathcal{H}(A,B,C)\cap D_2(\mathcal{F})=D_2(A,B,C),$$ we have 
\begin{equation}\label{d2+d3_eq_5}
|D_2(\mathcal{F},\mathcal{M})|=\frac{1}{(n-2)}\sum_{\{A,B,C\}\subseteq \mathcal{M}}|D_2(A,B,C)|.
\end{equation}
By the assumption that $|D_2(A,B,C)|\leq 21$ for all $\{A,B,C\}\subseteq \mathcal{M}$ and by equations (\ref{d2+d3_eq_4}) and (\ref{d2+d3_eq_5}), we get
\begin{eqnarray*}
|D_2(\mathcal{F},\mathcal{M})|+|D_3(\mathcal{F})\setminus \mathcal{M}|&\leq& \frac{2}{3}|D_2(\mathcal{F},\mathcal{M})|+\frac{27}{3(n-2)}{n \choose 3}\\
&=& {\frac{2}{3}\left(\frac{1}{(n-2)}\sum_{\{A,B,C\}\subseteq \mathcal{M}}|D_2(A,B,C)|\right)+\frac{27}{3(n-2)}{n \choose 3}}\\
&\leq & \frac{2}{3}\left(\frac{1}{(n-2)}\sum_{\{A,B,C\}\subseteq \mathcal{M}}21\right)+\frac{27}{3(n-2)}{n \choose 3}\\
&=& \left(\frac{2}{3}\right)\frac{21}{(n-2)} {n \choose 3}+\frac{27}{3(n-2)}{n \choose 3}\\
&=& \frac{23}{(n-2)}{n \choose 3}.
\end{eqnarray*}\hspace{\stretch{1}}$\square$

Let $\mathcal{M}_1$, $m$ and $\mathcal{M}_2$ be defined by Definition \ref{rem matching partition}. Also, define $n:=|\mathcal{M}_2|$.
\begin{prop}\label{prop bf-4}
Let $\mathcal{F}$ be a $3$-uniform, linear family and let $\mathcal{M}$ be a maximum matching of $\mathcal{F}$. If $\mathcal{E}_4$ is defined by Definition \ref{partition of family} then 
\begin{equation}
|\mathcal{E}_4|\leq
\begin{dcases} 
  \frac{23 n(n-1)}{6} , &\text{if\ } n\geq 3 \\
  8, &\text{if\quad} n=2\\
  0, &\text{if\quad $n=1$ \ or \ $n=0$.} 
 \end{dcases}
 \end{equation}
\end{prop}
\noindent{\sl Proof.} Let $n\geq 3$ and suppose that $|\mathcal{E}_4|> \frac{23}{6(n-2)}{n \choose 3}=\frac{23}{6}n(n-1)$. So, by Proposition \ref{prop_imp d2+d3}, there are edges $A$, $B$ and $C$ in $\mathcal{M}_2$ such that $|D_2(A,B,C)|>21$. But then there is an $\mathcal{M}_2$-augmenting set $\mathcal{W}$ in $\mathcal{F}$ by Proposition \ref{prop_c000} such that $\mathcal{W}\cap \mathcal{M}_2=\{A,B,C\}$ and $\mathcal{W}\setminus \mathcal{M}_2 \subset D_2(A,B,C)$. If edges in $\mathcal{W}\setminus \mathcal{M}_2$ do not intersect with any edge in $\mathcal{M}_1$ then $\mathcal{W}$ is an $\mathcal{M}$-augmenting set too. Thus, we have a contradiction to the fact that $\mathcal{M}$ is a maximum matching. So, there are edges in $\mathcal{W}$ that intersect with $X_{\mathcal{M}_1}$. By Definition \ref{augmenting set}, $|\mathcal{W}\setminus \mathcal{M}_2|\geq 4$. Let $B_1$, $B_2$, $B_3$ and $B_4$ be edges in $\mathcal{W}\setminus \mathcal{M}_2$. Note that if  $X_{\mathcal{M}_1}\cap B_i \neq \emptyset$ for some $i \in \{1,2,3,4\}$, then $B_i\in D_3(\mathcal{F},\mathcal{M})\cap \mathcal{E}_4$. Let $j:=|\{i \quad| \quad B_i \cap X_{\mathcal{M}_1}\neq \emptyset\}|$. By definition $0\leq j \leq 4$, so we need to consider cases for $j \in \{0,1,2,3,4\}$. In case $j=0$, the result is already established. One can easily construct an $\mathcal{M}$-augmenting set (similar to Proposition \ref{prop bf-3}) by considering $D_1(\mathcal{F})$ edges incident to $(\mathcal{M}_1)_{\mathcal{W}}$ edges in all cases for $j \in \{1,2,3,4\}$. Note that at most four edges in $\mathcal{M}_1$ can have non-empty intersection with $\cup_{i=1}^4B_i$. We leave details of construction of augmenting set for each case $j \in \{1,2,3,4\}$ to the readers. \newline
If $n=2$ then by Proposition \ref{prop_c00} and definition of $\mathcal{E}_4$, we have $|\mathcal{E}_4|\leq 8$. Also by Definition \ref{partition of family}, $\mathcal{E}_4$ is empty if $n<2$.
\hspace{\stretch{1}}$\square$ 

\medskip

We recall Definition \ref{rem matching partition} regarding partition of $\mathcal{M}$. In the proof of the following proposition, $m:=|\mathcal{M}_1|$, $\nu:=\nu(\mathcal{F})$ and $\Delta:=\Delta(\mathcal{F})$.
\begin{prop}\label{final_prop} Let $\mathcal{F}$ be a $3$-uniform, linear family such that $S_{\mathcal{F}}=\emptyset$, i.e., there is no vertex in $\mathcal{F}$ that is covered by all maximum matchings. If $\nu(\mathcal{F})=\nu$ then
 \begin{equation}\label{final_eq_1}
 |\mathcal{F}|\leq \frac{23}{6}{\nu}^2+7\nu.
 \end{equation}
 \end{prop}
\noindent{\sl Proof.} Let $\Delta:=\Delta(\mathcal{F})$. By Proposition \ref{degree_knu}, $S_{\mathcal{F}}=\emptyset$ implies that $\Delta\leq 3\nu$. By Definition \ref{partition of family} of $\mathcal{E}_i$'s, $|\mathcal{F}|\leq \sum_{i=1}^4{|\mathcal{E}_i|}+|\mathcal{M}_2|$. Proposition \ref{prop bf-1} implies that $|\mathcal{E}_1|\leq m\Delta$, Proposition \ref{prop bf-2} implies that $\mathcal{E}_2=\emptyset$, Proposition \ref{prop bf-3} implies that\newline $|\mathcal{E}_3|\leq (\nu-m)\min\{(2m+6),\Delta-1\}\leq (\nu-m)(2k+6)$ and by Proposition \ref{prop bf-4},
 $|\mathcal{E}_4|\leq \frac{23}{6}(\nu-m)(\nu-m-1)\leq \frac{23}{6}(\nu-m)^2$ for $\nu-m\geq 3$. Note that $|\mathcal{E}_4|\leq 8$ for $\nu-m\leq 2$. Also, $|\mathcal{M}_2|=\nu-m$. If $\nu-m\geq 3$, then\newline
\begin{eqnarray*}
|\mathcal{F}|&\leq & \sum_{i=1}^4{|\mathcal{E}_i|}+|\mathcal{M}_2|\\
             &\leq & m\Delta+(\nu-m)(2m+6)+\frac{23}{6}(\nu-m)^2+(\nu-m)\\
             &\leq & 3\nu m+(\nu-m)(2m+7)+ \frac{23}{6}(\nu-m)^2 \hspace{ 1 cm}{[\mbox{as $\Delta\leq 3\nu$}]}\\
             &= & m^2\left(-2+\frac{23}{6}\right)+ m\left(3\nu+2\nu-7-\frac{23}{3}\nu\right)+ \frac{23}{6}{\nu}^2+7\nu\\
             &=&  m^2\left(\frac{11}{6}\right)-m \left(\frac{8}{3}\nu+7\right)+\frac{23}{6}{\nu}^2+7\nu.\\
\end{eqnarray*}
The final expression above is a concave upward parabola in $m$ and hence the maximum value would occur at the extreme points, $m=0$ or $m=\nu-3\leq \nu$. It is easily checked that maximum occurs at $m=0$. Hence,
\begin{equation}\label{final_eq_2}
|\mathcal{F}|\leq \frac{23}{6}{\nu}^2+7\nu.
\end{equation}
If $\nu -m \leq 2$ then by Proposition \ref{prop bf-4}, $|\mathcal{E}_4|\leq 8$. Hence
\begin{eqnarray*}
|\mathcal{F}|&\leq & \sum_{i=1}^4{|\mathcal{E}_i|}+|\mathcal{M}_2|\\
             &\leq & m\Delta+(\nu-m)(\Delta -1)+8+(\nu-m) \mbox{ [as $|\mathcal{E}_3|\leq (\Delta-1)(\nu-m)$]}\\
             &= & \Delta \nu+8\\
             &\leq& 3\nu^2+8 \mbox{ [ as $\Delta \leq 3\nu$]}\\
             & \leq & \frac{23}{6}{\nu}^2+7\nu \mbox{ [ for $\nu \geq 2$]}.\end{eqnarray*}
For $\nu(\mathcal{F})=1$ and $\Delta(\mathcal{F})\leq 3\nu(\mathcal{F})=3$, use equation (\ref{trivial-bound}) to obtain $|\mathcal{F}|\leq 3\Delta-2\leq 7\leq \frac{23}{6}{\nu}^2+7\nu$.\hspace{\stretch{1}} $\square$\newline
\section{Proof of the main result-Theorem \ref{thm1}}
{\sl {\bf{Proof of Theorem \ref{thm1}}}}: Let $x \in X_{\mathcal{F}}$ be such that $|\mathcal{F}_x|=\Delta$. By Proposition \ref{degree_knu}, $x \in S_{\mathcal{F}}$ as $\Delta\geq\frac{23}{6}\nu(1+\frac{1}{\nu-1})> 3\nu$. Recall Definition \ref{nested-sequence}.  As $S_{\mathcal{F}}\neq \emptyset$, therefore there is a \textit{nested} sequence $\{y_1,\ldots,y_{k_1}\}\subseteq X_{\mathcal{F}}$. By Proposition \ref{sf-bound},\newline
\begin{equation}\label{final_thm_eq_1}
|\mathcal{F}|\leq k_1 \Delta+|\mathcal{F}_{k_1}|.
\end{equation}
Note that Proposition \ref{sf-bound} also implies that $\Delta(\mathcal{F}_{k_1})\leq 3 \nu(\mathcal{F}_{k_1})$. By the definition of $y_i$'s and repeated use of Remark \ref{SF-rem}, we get $\nu(\mathcal{F}_{k_1})=\nu-k_1$. Since $S_{\mathcal{F}_{k_1}}=\emptyset$, by Proposition \ref{final_prop} and equation (\ref{final_thm_eq_1}) we have\newline  
\begin{eqnarray*}
|\mathcal{F}|&\leq& k_1 \Delta+|\mathcal{F}_{k_1}|\\
             & \leq& k_1\Delta+ \frac{23}{6}(\nu-k_1)^2+7(\nu-k_1)\\
             & =& {k_1}^2(\frac{23}{6})+k_1(\Delta-\frac{23}{3}\nu-7)+\frac{23}{6}{\nu}^2+7\nu.
\end{eqnarray*}  
Let $f(k_1):={k_1}^2(\frac{23}{6})+k_1(\Delta-\frac{23}{3}\nu-7)+\frac{23}{6}{\nu}^2+7\nu$ for $1\leq k_1\leq \nu$. Note that $k_1\geq 1$ because $S_{\mathcal{F}}\neq \emptyset$. Clearly $f(k_1)$ is a concave upward parabola as $\frac{d^2f(k_1)}{d{k_1}^2}>0$. Hence the maximum of $f(k_1)$ occurs at the extreme points $k_1=1$ or $k_1=\nu$. As 
$f(1)=\frac{23}{6}+\Delta-\frac{23}{3}\nu-7+\frac{23}{6}{\nu}^2+7\nu=\frac{23}{6}{\nu}^2+\Delta-\frac{2\nu}{3}-\frac{19}{6}\leq \frac{23}{6}{\nu}^2+\Delta$ and $f(\nu)=\Delta \nu$. Thus, $|\mathcal{F}|\leq \max\{\frac{23}{6}{\nu}^2+\Delta,\Delta \nu \}$. Since $\Delta \nu\geq \frac{23}{6}{\nu}^2+\Delta$ if and only if $\Delta \geq \frac{23}{6}\frac{{\nu}^2}{(\nu-1)}$. Therefore for $\Delta \geq \frac{23}{6}\frac{{\nu}^2}{(\nu-1)}$,  \newline
\begin{displaymath}
|\mathcal{F}|\leq \Delta \nu. 
\end{displaymath}\hspace{\stretch{1}} $\square$

Recall by Remark \ref{delta-nu}, for any positive integers $\Delta$ and $\nu$ there exists a $3$-uniform, linear family $\mathcal{F}$ with $\Delta(\mathcal{F})=\Delta$, $\nu(\mathcal{F})=\nu$ such that $|\mathcal{F}|=\Delta \nu$. Thus, an extremal family achieves the bound on the size in Theorem \ref{thm1}.
\section{Acknowledgements}
The author would like to thank Dr. Nishali Mehta and Dr. Naushad Puliyambalath for their valuable comments. This article is part of author's doctoral research that was guided by Prof. Akos Seress. The author is indebted to his advisor for suggesting the problem, sharing critical insights and steering the course of the research.      

\end{document}